\pdfoutput=1\relax
\documentclass[11pt, english]{article}
\usepackage{meta/preamble}
\usepackage{mathtools}
\usepackage{tikz-cd}

\makeatletter
\@ifundefined{example}{\theoremstyle{definition}\newtheorem{example}[theorem]{Example}}{}
\@ifundefined{corollary}{\theoremstyle{plain}\newtheorem{corollary}[theorem]{Corollary}}{}
\makeatother
\theoremstyle{plain}
\newtheorem{introthm}{Theorem}

\providecommand{\cofib}{\operatorname{cofib}}
\providecommand{\Csc}{\Cond(\An)_*^{\ge2}}
\providecommand{\Csol}{\Cond(\An)_{*,\mathrm{sol}}^{\ge2}}

\title{Unstable Solidification}
\date{\today}
\author{Shachar Carmeli, Amos Kaminski\thanks{Department of Mathematics, Weizmann Institute of Science \\shachar.carmeli@weizmann.ac.il \\
amos.kaminski@weizmann.ac.il}}
\begin{document}
\maketitle
\begin{abstract}
We construct a solidification functor $\LH$ on condensed anima by
localizing at the maps inverted by solid integral chains. A pointed
connected nilpotent condensed anima is $\Hsol$-local precisely when its
fundamental group admits a finite central series with solid abelian
factors and its higher homotopy groups are solid. For simply connected
condensed anima, we construct a natural strongly convergent spectral sequence
\[
E^2_{p,q}=L_ps(\pi_qX)\Longrightarrow\pi_{p+q}\LH X
\]
 where $L_ps$ denotes
the $p$-th derived functor of the abelian solidification.
\end{abstract}
\begin{center}
    \includegraphics[width=0.58\textwidth]{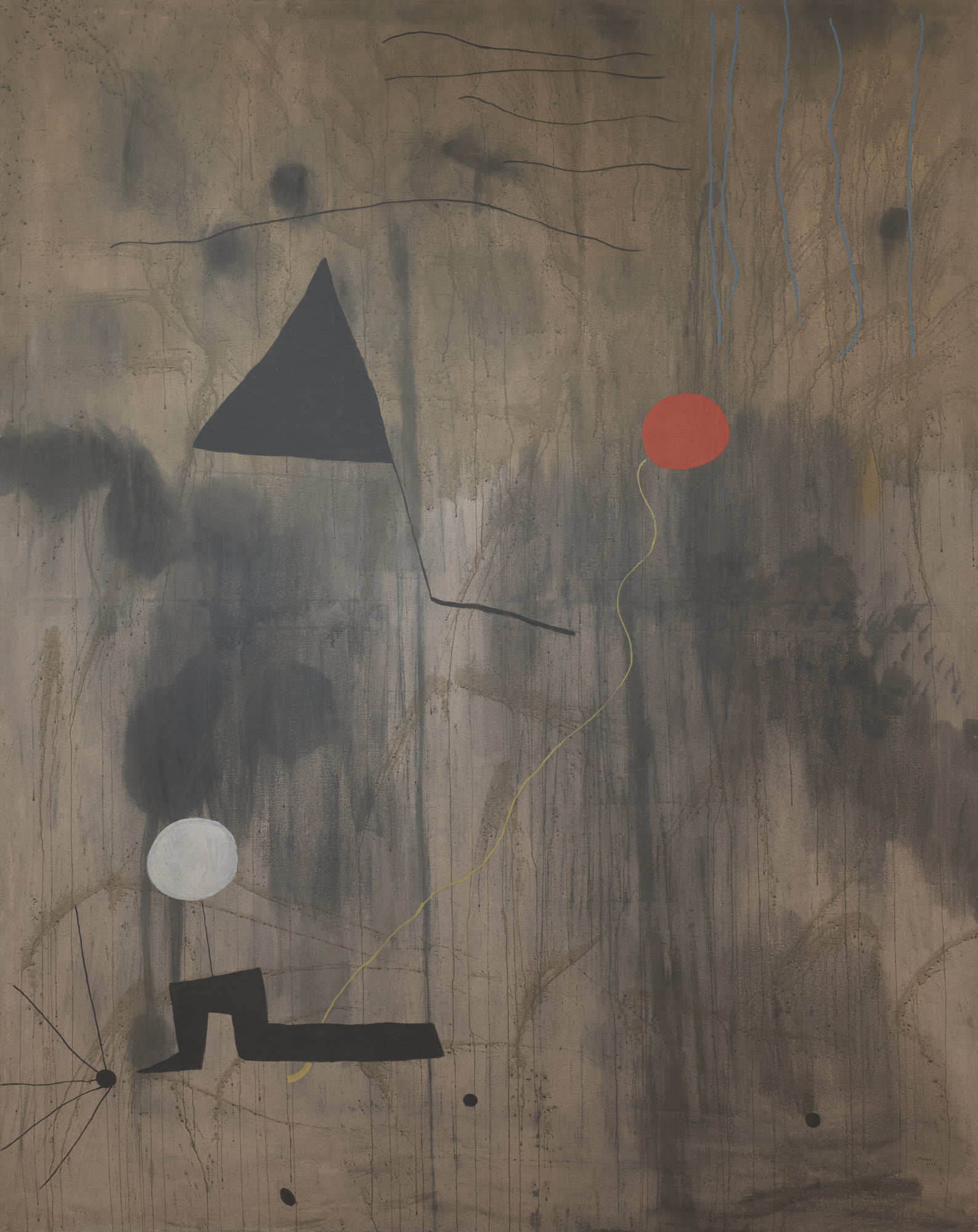}
    
    \smallskip
    {\footnotesize
    Joan Miró, \emph{La Naissance du monde} (1925).\\
    The Museum of Modern Art, New York.
    }
\end{center}

\section*{Introduction}
Solidity, introduced by Clausen and Scholze, is a non-archimedean notion of completeness for condensed abelian groups, and its theory has been extended to condensed spectra. 
The solidification functor, is the left
adjoint $s$
of the inclusion of solid abelian groups into condensed abelian
groups. It has a well-behaved left derived functor $\Ls$ and a stable
counterpart for condensed spectra \cite{Aoki}. The purpose of this paper is to
construct its unstable counterpart for condensed anima, and to compute it. 
We will follow Bousfield's classical construction
of homological localization: a generalized homology theory $h_*$
determines a localization of spaces at the maps inducing
isomorphisms on $h_*$
\cite[Theorem~3.2]{Bousfield1975}.
Its local spaces are characterized by the requirement that
mapping into them turns every $h_*$-equivalence into an equivalence
of mapping spaces.  We apply the same principle to condensed anima,
using solid integral chains to detect equivalences  (Definition~\ref{def:Hlocal}).
Namely, we define $\LH$ to be the localization of $\Cond(\An)$ at the
maps inverted by the solid integral chain functor
$\Hsol=\Ls\Z[-]$.

The effect of this localization on homotopy groups cannot, in
general, be described by applying $s$ degreewise. It fails for two reasons.
First, $s$ is not exact, so its derived functors  
$L_ps$ must intervene. For example, for
$\mathbb T=\underline{\mathbb R/\Z}$ one has $s\mathbb T=0$ but
$\Ls\mathbb T\simeq\Z[1]$, and the solidification of $K(\mathbb T,n)$ turns
out to be $K(\Z,n+1)$. Second, a homotopy type is more than its homotopy
groups: these are assembled by the Postnikov $k$-invariants, and a
solidification has to carry these attaching maps along.

However, at least in the simply connected case, the $\Hsol$-local condensed anima can be characterized by their homotopy groups.
 Write
$\Cond(\An)_*^{\ge2}$ for the pointed simply connected condensed anima and
$\Cond(\An)_{*,\mathrm{sol}}^{\ge2}$ for those whose homotopy groups are
solid.

\begin{introthm}[Theorem~\ref{thm:reflection}]\label{introthm:reflection}
A pointed simply connected condensed anima is $\Hsol$-local if and only if all its
homotopy groups are solid. Moreover, $\LH$ preserves simply connected
condensed anima. Hence the restriction of $\LH$ to simply connected
condensed anima is their reflection into the
condensed anima with solid homotopy groups.
\end{introthm}

Having identified the $\Hsol$-local simply connected condensed anima, we turn to the effect of
localization on homotopy groups. Derived solidification computes
the Eilenberg--Mac Lane case and supplies the terms of a spectral
sequence for general simply connected condensed anima.

\begin{introthm}[Proposition~\ref{prop:EM} and Theorem~\ref{thm:SS}]\label{introthm:SS}
For a $1$-connective  condensed spectrum $E$ one has
$\LH \Omega^\infty E\simeq\Omega^\infty\Lsp E$, where $\Lsp$ is the solidification
of spectra \ref{prop:stable}. In particular, for a condensed abelian group $A$ and $n\ge2$,
\[
 \LH K(A,n)\simeq\Omega^\infty\UZ\Ls(A[n]),
 \qquad
 \pi_{n+p}\LH K(A,n)\simeq L_ps(A)\quad(p\ge0)
\]
where $\UZ$ denotes the forgetful functor from condensed $H\Z$ modules into spectra.
For every $X\in\Cond(\An)_*^{\ge2}$ there is a natural strongly convergent
spectral sequence
\[
 E^2_{p,q}=L_p s(\pi_qX)\Longrightarrow\pi_{p+q}\LH X,
 \qquad p\ge0,\quad q\ge2,
\]
with $d_r:(p,q)\mapsto(p-r,q+r-1)$ and with a finite filtration in each total
degree.
\end{introthm}
We construct the spectral sequence by solidifying the Whitehead tower.

Beyond the simply connected case one has to keep track of $\pi_1$: both $B\pi_1X$ and the action of $\pi_1X$ on the higher homotopy groups interact with solidification. In the nilpotent case we nevertheless obtain a complete characterization: 

\begin{introthm}[Corollary~\ref{cor:nilpotent}]\label{introthm:nilpotent}
Let $X$ be a pointed connected nilpotent condensed anima. Then $X$ is
$\Hsol$-local if and only if $\pi_1(X)$ admits a finite central series with solid factors \footnote{Recall that a finite central series for $G$ is a filtration
$1=G_0\subseteq\cdots\subseteq G_c=G$ by normal subgroups such that
each $G_i/G_{i-1}$ lies in the center of $G/G_{i-1}$.} \footnote{We require only that $G$ admit one central series with solid
factors; other central series of $G$ need not have this property.} and $\pi_nX$ is solid for every $n\ge2$. Moreover $\LH X$ is
again connected and nilpotent.
    
\end{introthm}

The condition on $\pi_1X$ cannot be weakened to solidity of its abelianization: a Heisenberg
group gives a counterexample (Example~\ref{ex:heisenberg}). Since
discrete abelian groups are solid, every nilpotent constant condensed anima is
$\Hsol$-local. Consequently a connected finite CW complex with nilpotent homotopy type,
viewed as a condensed set (i.e., static condensed anima), localizes to its homotopy type viewed as a \emph{discrete but non-static} condensed anima. For instance, 
\[
\LH\mathbb T\simeq K(\Z,1),
\]
where the LHS is static but non-discrete while the RHS is discrete but non-static.

Finally, one should note that pointed simply connected solid anima are not closed under colimits, in contrast with the case of solid spectra.
\begin{example}[pointed simply connected solid anima are not closed under colimit-Proposition~\ref{prop:noadjoints}]
    For $P=\prod_{\mathbb N}\Z$, the wedge $Q:=K(P,2)\vee K(P,2)$ has
$\pi_3Q=P\otimes P$, which is not solid; hence the inclusion of pointed, solid, simply
connected condensed anima into arbitrary pointed simply connected condensed anima has no right adjoint,
although its stable analogue does. The coproduct in the solid category has
$\pi_3Q=P\otimes_\square P$ instead (Example~\ref{ex:wedge}).
\end{example}

\subsection*{Acknowledgments}
We thank Dustin Clausen for suggesting the idea that led to this work.

\subsection*{Funding} 
This Research is partially supported by ISF Beresheet grant 4093/25, BSF grant 2024766, Minerva grant 715294, and the Azrieli foundation through an Early Career Researcher Fellowship.
\subsection*{Acknowledgment of AI assistance}
The mathematical ideas and research direction are due to the authors. OpenAI's
ChatGPT/Codex and Anthropic's Claude were used under their direction to assist
with drafting and revision, literature searches, and the development and
checking of arguments. This assistance included generating portions of the
technical proofs. The authors take full responsibility for the mathematical
content.

\section{Background}
\label{sec:solid-background}
\subsection{Solid abelian groups}
Following Clausen and Scholze, we begin by recalling some basic definitions and facts regarding solid abelian groups and their derived category. 
\begin{definition}[Solid abelian groups {\cite[Definition 5.1]{CS}}]
For a profinite set $S=\varprojlim_i S_i$, with $S_i$ finite, set
\[
 \Z[S]^\square:=\varprojlim_i\Z[S_i].
\]
A condensed abelian group $A$ is \emph{solid} if the canonical map
$\Z[S]\to\Z[S]^\square$ induces a bijection
\[
 \operatorname{Hom}_{\Cond(\Ab)}(\Z[S]^\square,A)
 \xrightarrow{\;\sim\;}
 \operatorname{Hom}_{\Cond(\Ab)}(\Z[S],A)
\]
for every profinite set $S$. Write $\Solid(\Ab)$ for the full subcategory of
solid condensed abelian groups.
\end{definition}
Recall the following structural results for solid abelian groups due to Clausen and Scholze:
\begin{theorem}[Clausen--Scholze {\cite[Theorem 5.8]{CS}}]\label{thm:solid-abelian}
The following holds for $\Solid(\Ab)$.
\begin{enumerate}
\item
$\Solid(\Ab)$ is an abelian subcategory of $\Cond(\Ab)$ closed under all limits,
colimits, and extensions.
\item
Products of copies of $\Z$ form a family of compact projective generators of
$\Solid(\Ab)$.
\item
The inclusion $\Solid(\Ab) \hookrightarrow \Cond(\Ab) $ admits a left adjoint
\[
 s:\Cond(\Ab)\longrightarrow\Solid(\Ab)
\]
which satisfies 
\[
s(\Z[S])\simeq\Z[S]^\square.
\]
\end{enumerate}
\end{theorem}

In particular $\Z$ is solid, and hence so is every discrete abelian group, being a cokernel of a map between direct sums of copies of $\Z$. Since $s$ is a left adjoint, it is right exact; its failure to be exact is measured by its left derived functor:

\begin{theorem}[Derived solidification {\cite[Theorem 5.8]{CS}}]\label{thm:solid-derived}
The functor
\[
 D(\Solid(\Ab))\longrightarrow D(\Cond(\Ab))
\]
is fully faithful. Its essential image consists of the complexes whose homology
groups are solid. Equivalently, $D(\Solid(\Ab))$ consists of  the complexes $C$ for which the map
\[
R\operatorname{Hom}(\Z[S]^\square,C)\to R\operatorname{Hom}(\Z[S],C)
\]
is an
isomorphism for all profinite $S$.
The inclusion
\[
D(\Solid(\Ab)) \hookrightarrow D(\Cond(\Ab)) 
\]
admits a left adjoint
\[
 \Ls:D(\Cond(\Ab))\longrightarrow D(\Solid(\Ab))
\]
preserving connective complexes
which is the left derived functor of $s$.
\end{theorem}

\begin{remark}
We use homological indexing and write
\[
 L_p s(A):=H_p(\Ls(A[0])).
\]
Thus $L_0s(A)=sA$, and $L_ps(A)=0$ for $p<0$. If $A$ is solid, then
$\Ls(A[0])\simeq A[0]$, so $L_ps(A)=0$ for $p>0$.
\end{remark}

The inclusion in Theorem~\ref{thm:solid-derived} is $t$-exact, so
$D(\Solid(\Ab))$ inherits a $t$-structure from $D(\Cond(\Ab))$.

We next recall how solidification interacts with tensor products.
\begin{theorem}[Solid tensor products {\cite[Theorem 6.2 and Proposition 6.3]{CS}}]\label{thm:solid-tensor}
There are symmetric monoidal structures on $\Solid(\Ab)$ and $D(\Solid(\Ab))$
for which $s$ and $\Ls$ are symmetric monoidal, that is, their tensor products satisfy
\[
 A\otimes_\square B=s(A\otimes B),
 \qquad
 M\sotimes N=\Ls(M\otimes^{\mathbb L}N).
\]
Moreover,
\[
 \left(\prod_I\Z\right)\sotimes\left(\prod_J\Z\right)
 \simeq\prod_{I\times J}\Z.
\]
\end{theorem}

By adjunction, the monoidality of $\Ls$ implies that solid complexes are closed under internal mapping complexes out of arbitrary condensed complexes, namely:
\begin{lemma}\label{lem:internal-hom-solid}
If $M\in D(\Cond(\Ab))$ is solid and $A\in D(\Cond(\Ab))$ is arbitrary, then the
internal mapping complex $R\underline{\operatorname{Hom}}(A,M)$ is solid. In
particular, if $M$ is a solid abelian group and $A$ is a condensed abelian group,
then $\underline{\operatorname{Hom}}(A,M)$ is solid.
\end{lemma}

\begin{proof}
The monoidality of $\Ls$ gives isomorphisms

\[
\Ls(\Z[S]^\square\otimes^{\mathbb L}A)
\simeq\Z[S]^\square\sotimes\Ls A\simeq
\Ls\Z[S]\sotimes\Ls A 
\simeq\Ls(\Z[S]\otimes^{\mathbb L}A)
\]
for every profinite set $S$. Since $M$ is solid,
$R\operatorname{Hom}(B,M)\simeq R\operatorname{Hom}(\Ls B,M)$ for every $B$. Applying this to $B = \Z[S]\otimes^{\mathbb L}A$ we obtain:
\begin{align*}
 R\operatorname{Hom}(\Z[S]^\square,R\underline{\operatorname{Hom}}(A,M))
 &\simeq R\operatorname{Hom}(\Z[S]^\square\otimes^{\mathbb L}A,M)\\
 &\simeq R\operatorname{Hom}(\Z[S]\otimes^{\mathbb L}A,M)
 \simeq R\operatorname{Hom}(\Z[S],R\underline{\operatorname{Hom}}(A,M)).
\end{align*}
The result now follows from  Theorem~\ref{thm:solid-derived}. The second assertion follows from the first by
taking $H_0$.
\end{proof}
\subsection{Solid spectra}

Similarly to the derived condensed abelian groups, one can define a notion of solid condensed spectra. We will compare unstable solidification with its stable counterpart in Section~\ref{sec:stable}, so we recall it here.
 The characterization of solid complexes by their homology
suggests the following equivalent definitions for spectra.
\begin{definition}[Solid spectra
{\cite[Section 2.4]{Aoki}}]
\label{def:solid-spectra}
For a profinite set $S=\varprojlim_i S_i$
set
\[
\mathbb S[S]:=\Sigma^\infty_+S,
\qquad
\mathbb S[S]^\square:=\varprojlim_i\mathbb S[S_i].
\]
One has a canonical map
$\mathbb S[S]\to\mathbb S[S]^\square$.

The following conditions on a condensed spectrum $E$ are equivalent:
\begin{enumerate}
\item
The condensed abelian group $\pi_nE$ is solid for every $n\in\Z$.

\item
For every profinite set $S$, the canonical map induces an
equivalence of mapping spectra
\[
\operatorname{map}_{\Cond(\Sp)}(\mathbb S[S]^\square,E)
\xrightarrow{\;\sim\;}
\operatorname{map}_{\Cond(\Sp)}(\mathbb S[S],E).
\]

\item
The spectrum $E$ belongs to the smallest full stable subcategory
of $\Cond(\Sp)$ closed under colimits and containing all products
$\prod_I\mathbb S$ of the condensed sphere spectrum.
\end{enumerate}
A condensed spectrum satisfying these conditions is called
\emph{solid}. We write $\Solid(\Sp)\subseteq\Cond(\Sp)$ for
the full subcategory of solid spectra.

\end{definition}
\begin{proof}
The equivalence of (1) and (3) follows from
\cite[Theorem~2.19 and Definition~2.21]{Aoki}, using
$E\simeq\varinjlim_m\tau_{\ge-m}E$ for unbounded spectra.
By \cite[Theorem~5.4]{CS} $C(S,\Z)$ is a free abelian groups, hence we can choose a basis $C(S,\Z)\cong\bigoplus_I\Z$.
Remark that now the canonical map \[
f:\bigoplus_{b\in I}\mathbb S\longrightarrow\varinjlim_i\mathbb S^{S_i}.
\] is an equivalence since it's an equivalence on homotopy groups:
For every $k\in\Z$, we have
\[
\begin{aligned}
\pi_k\!\left(\varinjlim_i\mathbb S^{S_i}\right)
&\cong \varinjlim_i(\pi_k\mathbb S)^{S_i}
 \cong \varinjlim_i\bigl(\Z^{S_i}\otimes_\Z\pi_k\mathbb S\bigr) \\
&\cong \left(\varinjlim_i\Z^{S_i}\right)\otimes_\Z\pi_k\mathbb S
 \cong C(S,\Z)\otimes_\Z\pi_k\mathbb S \\
&\cong \left(\bigoplus_{b\in I}\Z\right)\otimes_\Z\pi_k\mathbb S
 \cong \bigoplus_{b\in I}\pi_k\mathbb S \\
&\cong \pi_k\!\left(\bigoplus_{b\in I}\mathbb S\right).
\end{aligned}
\]

Dualizing gives
\[
\mathbb S[S]^\square
=\varprojlim_i\mathbb S[S_i]
\simeq\underline{\operatorname{Hom}}\!
  \left(\varinjlim_i\mathbb S^{S_i},\mathbb S\right)
\simeq\underline{\operatorname{Hom}}\!
  \left(\bigoplus_I\mathbb S,\mathbb S\right)
\simeq\prod_I\mathbb S.
\].

Now remark that
\[
H\Z\otimes_{\mathbb S}\mathbb S[S]^\square
\simeq H\Z\otimes_{\mathbb S}\prod_I\mathbb S
\xrightarrow{\;\sim\;}\prod_I H\Z
\simeq H(\Z[S]^\square),
\]
where the middle map is an equivalence by \cite[Lemma~2.20]{Aoki}.
Hence the map
$\mathbb S[S]\to\mathbb S[S]^\square$ induces an equivalence
of mapping spectra into every shift of $HA$ with $A$ solid.
Postnikov induction and passage to the limit give the same conclusion
for every spectrum with solid homotopy groups.
This proves (1) implies 2.

Conversely, consider the full subcategory of spectra satisfying~(2). It is
reflective by \cite[\S5.5.4]{HTT}. Each $\mathbb S[S]^\square$
satisfies~(1), hence~(2), and condition~(2) says precisely that
$\mathbb S[S]\to\mathbb S[S]^\square$ is its reflection.
Since the shifts of $\mathbb S[S]$ for $S$ extremally disconnected
generate $\Cond(\Sp)$ under colimits, their reflections generate
this subcategory under colimits.

Now spectra satisfying~(1) are closed under coproducts and cofibers
by \cite[Theorem~5.8]{CS} and the homotopy exact sequence, hence
under all colimits. By $(1)\Rightarrow(2)$, these colimits already
satisfy~(2). Thus the spectra satisfying~(1) contain the generators
and are closed under colimits in this subcategory, proving
$(2)\Rightarrow(1)$.
\end{proof}
We now collect the formal properties of $\Solid(\Sp)$ that will be used below.
\begin{proposition}[cf {\cite[Theorem 2.19]{Aoki}}]\label{prop:stable}
The category $\Solid(\Sp)$ is a stable subcategory of $\Cond(\Sp)$ whose inclusion
preserves limits and colimits, and the standard $t$-structure restricts to it with
heart $\Solid(\Ab)$. The inclusion has an exact, connective-preserving left adjoint
$\Lsp$ and a right adjoint $R_{sp}$. 
\end{proposition}
\begin{proof}

The locality characterization exhibits $\Solid(\Sp)$ as an accessible
stable localization of $\Cond(\Sp)$, while the generation characterization
shows that its inclusion preserves colimits. The inclusion therefore
preserves limits and colimits and admits both adjoints by the adjoint
functor theorem. Its left adjoint $\Lsp$ is exact by stability.

The homotopy-group characterization shows that truncations preserve
solidity, so the standard $t$-structure restricts, with heart
$\Solid(\Ab)$. Since the inclusion is $t$-exact, its left adjoint
preserves connective spectra.
\end{proof}

The closure under internal Homs established for solid complexes also holds for solid spectra.
\begin{lemma}\label{solidhom}
    If $Y\in \Cond(\Sp)$ is solid and $W\in \Cond(\Sp)$, then the
internal mapping spectrum $\underline{\operatorname{Hom}}(W,Y)$ is solid.
\end{lemma}
\begin{proof}
Since solid spectra are closed under limits and shifts, and
$\underline{\operatorname{Hom}}(-,Y)$ takes colimits to limits,
it suffices to treat the generators $W=\mathbb S[S]$, where
$S$ is profinite,
which then is solid by
\cite[Theorem~2.19(2)]{Aoki}.
\end{proof}
\section{Unstable solidification}
\label{sec:unstable}
Solid integral chains 
\[
\Hsol:= \Ls\Z[-]:\Cond(\An)\longrightarrow D(\Solid(\Ab)).
\]determine a class of equivalences on
condensed anima. We construct the associated localization:
\begin{definition}[Homological solidification]\label{def:Hlocal}
The \emph{homological solidification} 
\[
\LH: \Cond(\An) \to \Cond(\An)
\] 
is the localization of $\Cond(\An)$
at the maps inverted by the solid integral chain functor $\Hsol$.
Such maps are called \emph{$\Hsol$-equivalences}. A condensed anima $Y\in \Cond(\An)$ is
\emph{$\Hsol$-local} if $\Map(V,Y)\to\Map(U,Y)$ is an equivalence for every
$\Hsol$-equivalence $U\to V$. Thus $\LH$ reflects into the full subcategory
$\Cond(\An)_{\Hsol}$ of $\Hsol$-local condensed anima.
\end{definition}
We first check that this localization exists and is compatible with finite products.
\begin{proposition}\label{prop:H-exists}
The functor $\LH$ exists and preserves finite products.
\end{proposition}

\begin{proof}
Both integral chains and derived solidification are left adjoints,
so
\[
\Hsol=\Ls\Z[-]:\Cond(\An)\longrightarrow D(\Solid(\Ab))
\]
preserves colimits. By \cite[Proposition~5.5.4.16]{HTT}, the
class of $\Hsol$-equivalences is strongly saturated and generated
by a set of morphisms. Hence
\cite[Proposition~5.5.4.15]{HTT} gives the accessible localization
$\LH$.

To prove preservation of products, first observe that the integral
chain functor $\Z[-]$ is symmetric monoidal from cartesian products
to derived tensor products.
Derived solidification $\Ls$ is also symmetric monoidal, by
Theorem~\ref{thm:solid-tensor}. Consequently,
\[
\Hsol(X\times Y)
\simeq
\Ls\bigl(\Z[X]\otimes^{\mathbb L}\Z[Y]\bigr)
\simeq
\Hsol(X)\sotimes\Hsol(Y).
\]

Let $\eta_X:X\to\LH X$ and $\eta_Y:Y\to\LH Y$ be the localization
maps. The displayed equivalence identifies
$\Hsol(\eta_X\times\eta_Y)$ with the tensor product of the
equivalences $\Hsol(\eta_X)$ and $\Hsol(\eta_Y)$. Thus
\[
X\times Y\longrightarrow\LH X\times\LH Y
\]
is an $\Hsol$-equivalence. Its target is $\Hsol$-local, since mapping
into a product gives the product of the mapping spaces.
It therefore exhibits $\LH X\times\LH Y$ as the localization
of $X\times Y$, proving
\[
\LH(X\times Y)\simeq\LH X\times\LH Y.
\]
Finally, the terminal object is $\Hsol$-local because every mapping
space into it is contractible. Hence $\LH$ preserves finite
products.

\end{proof}
We will also need a pointed version of the
localization, defined using reduced chains.
\begin{definition}\label{def:pointed-H-equivalences}
For a pointed condensed anima $X$ with base point $x$, set
\[
\widetilde\Z[X]
:=\cofib\bigl(\Z\xrightarrow{x_*}\Z[X]\bigr),
\qquad
\widetilde\Hsol(X):=\Ls\widetilde\Z[X].
\]
A pointed map $f$ is a \emph{pointed $\Hsol$-equivalence} if
$\widetilde\Hsol(f)$ is an equivalence. A pointed condensed anima $Y$
is \emph{pointed $\Hsol$-local} if
$\Map_*(B,Y)\to\Map_*(A,Y)$ is an equivalence for every
pointed $\Hsol$-equivalence $A\to B$.
\end{definition}
The pointed version introduces nothing new: it is computed on underlying condensed anima.
\begin{proposition}\label{prop:pointed-localization}
Pointed $\Hsol$-equivalences and pointed $\Hsol$-local condensed anima
are detected by the forgetful functor $U\colon \Cond(\An)_*\to\Cond(\An)$. The pointed localization with respect to the \emph{pointed $\Hsol$-equivalences} exists  and is given by
\[
\LH^{\mathrm{pt}}(X,x)
=\bigl(\LH(UX),\eta_{UX}\circ x\bigr).
\]
In particular, $U\LH^{\mathrm{pt}}\simeq\LH U$.
\end{proposition}

\begin{proof}
The augmentation $\Z[X]\to\Z$ splits the base-point inclusion,
giving a natural decomposition
\[
\Hsol(X)\simeq\Z\oplus\widetilde\Hsol(X).
\]
Thus pointed $\Hsol$-equivalences are detected by $U$.

If $Y$ is $\Hsol$-local, then $Y$ is pointed $\Hsol$-local by taking fibers over
its base point in
\[
\Map_*(A,Y)\longrightarrow\Map(UA,UY)\longrightarrow\Map(*,UY).
\]
Conversely, for $V_+=V\amalg *$, pointed by the second summand,
\[
\widetilde\Hsol(V_+)\simeq\Hsol(V),
\qquad
\Map_*(V_+,Y)\simeq\Map(V,UY).
\]
Testing pointed $\Hsol$-locality on maps $A_+\to B_+$ therefore shows
that $Y$ is $\Hsol$-local.

Finally, the pointed unit
$(X,x)\to(\LH(UX),\eta_{UX}\circ x)$ is a pointed
$\Hsol$-equivalence with pointed $\Hsol$-local target, so it exhibits
the asserted localization.
\end{proof}

We henceforth write $\LH$ also for $\LH^{\mathrm{pt}}$.
The chain adjunction supplies our first examples of
$\Hsol$-local condensed anima. For $M\in D(\Solid(\Ab))$, it gives 
\begin{equation}\label{eq:chainadj-general}
 \Map(X,\Omega^\infty M)\simeq\Map(\Hsol(X),M),
 \qquad
 \Map_*(X,\Omega^\infty M)\simeq\Map(\widetilde\Hsol(X),M),
\end{equation}
the mapping spaces on the right being taken in $D(\Solid(\Ab))$.
The first equivalence shows that the underlying anima $\Omega^\infty M$ is
$\Hsol$-local. In particular, for $A$ solid and $n\ge0$, the condensed anima $K(A,n)$ is
$\Hsol$-local\label{EM_is_solid}.
For maps between simply connected condensed anima, $\Hsol$-equivalences are detected by the
condensed anima $K(A,n)$ with $A$ solid and $n\ge2$. Indeed, \eqref{eq:chainadj-general} reads
\begin{equation}\label{eq:chainadj}
 \Map_*(X,K(A,n))
 \simeq\Map_{D(\Solid(\Ab))}(\widetilde\Hsol(X),A[n]).
\end{equation}
If the cofiber of a map of such chains is nonzero, it is $2$-connective, and its
projection onto its first nonzero homology group detects it. Thus these tests detect
precisely the maps inverted by $\LH$.\\ 
Our main computational tool is the following analogue of the Bousfield--Kan fiber lemma.
\begin{proposition}[Fiber lemma]\label{prop:fiber}
If $F\to E\to B$ is a pointed fiber sequence of condensed anima and $B$ is simply connected, then
$\LH F\to\LH E\to\LH B$ is a fiber sequence.
\end{proposition}

\begin{proof}
This is the Bousfield--Kan fiber lemma \cite[Theorem 4.3]{BEM}, applied to the
functor $\Hsol$. Its hypotheses hold: $\Cond(\An)$ is an $\infty$-topos, $\Hsol$ is a symmetric monoidal left adjoint with connective values, and
$D(\Solid(\Ab))$ is presentably symmetric monoidal and stable, with a left-separated
$t$-structure. This $t$-structure is compatible with the tensor product, since $\Ls$
preserves connective complexes. The theorem requires both $B$ and
$\LH B$ to be simply connected, and the second condition is supplied by
Theorem~\ref{thm:reflection} below.
\end{proof}

\section{Characterizing the $\Hsol$-local condensed anima}
\label{sec:local-objects}

We now identify the simply connected $\Hsol$-local condensed anima.

Recall that $\Csc$ is closed under colimits in $\Cond(\An)_*$.
It is presentable, being the fiber of the colimit-preserving functor $\tau_{\le1}$
over the zero object. The inclusion has as right adjoint
the simply connected cover
\[
 \tau_{\ge 2}Z=\fib(Z\to\tau_{\le1}Z).
\]
So limits in $\Csc$ are the simply connected covers of the limits in
$\Cond(\An)_*$, and $\tau_{\ge 2}$ changes no $\pi_n$ with $n\ge2$.

\begin{theorem}\label{thm:reflection}\label{thm:homological}
A pointed simply connected condensed anima is $\Hsol$-local if and only if all its
homotopy groups are solid. Moreover, $\LH$ preserves simply connected
condensed anima. Hence the restriction of $\LH$ to simply connected
condensed anima is their reflection into the
condensed anima with solid homotopy groups.
\end{theorem}

\begin{proof}
We first show that simply connected condensed anima with solid homotopy groups are $\Hsol$-local.
Let $Z\in\Csol$. We prove by induction on $n\ge1$ that
$\tau_{\le n}Z$ is $\Hsol$-local.

For $n=1$, simple connectedness gives $\tau_{\le1}Z\simeq *$,
which is $\Hsol$-local. Now let $n\ge2$ and assume that
$\tau_{\le n-1}Z$ is $\Hsol$-local. The Postnikov extension
$\tau_{\le n}Z\to\tau_{\le n-1}Z$ is a gerbe whose band is
constant with value $\pi_nZ$, since its base is simply connected.
It is therefore classified by a $k$-invariant
\[
\kappa_n:\tau_{\le n-1}Z\longrightarrow K(\pi_nZ,n+1),
\]
with
\[
\tau_{\le n}Z
\simeq
\tau_{\le n-1}Z\times_{K(\pi_nZ,n+1)}*
\]
\cite[Section~7.2.2.26]{HTT}.
Since $\pi_nZ$ is solid, $K(\pi_nZ,n+1)$ is $\Hsol$-local\ref{EM_is_solid}.
The induction hypothesis and closure of $\Hsol$-local condensed anima under
pullbacks imply that $\tau_{\le n}Z$ is $\Hsol$-local.

Thus every Postnikov truncation of $Z$ is $\Hsol$-local. Postnikov
completeness of $Cond(An)$ and closure of $\Hsol$-local condensed anima under limits now give
that $Z\simeq\liml_n\tau_{\le n}Z$ is $\Hsol$-local.

Now we turn to show the converse namely that {simply connected $\Hsol$-local condensed anima have solid homotopy groups.}

Let $X$ be a pointed simply connected condensed anima, and let
$\eta\colon X\to Y$ be its reflection into the full subcategory of
pointed simply connected condensed anima with solid homotopy groups.
We claim that $\eta$ is an $\Hsol$-equivalence.

Fix $M\in D(\Solid(\Ab))$. The simply connected cover
$c_2\Omega^\infty M$ has solid homotopy groups: in degree $n\ge2$,
its homotopy group is $H_n(M)$. Since $X$ and $Y$ are simply
connected, pointed maps from either into $\Omega^\infty M$ factor
through its simply connected cover. Thus the reduced-chain adjunction
gives a commutative square
\[
\begin{array}{ccc}
\Map(\widetilde\Hsol(Y),M) &\longrightarrow&
\Map(\widetilde\Hsol(X),M)\\
\big\downarrow{\scriptstyle\simeq} && \big\downarrow{\scriptstyle\simeq}\\
\Map_*(Y,c_2\Omega^\infty M) &\xrightarrow{\;\eta^*\;}&
\Map_*(X,c_2\Omega^\infty M).
\end{array}
\]
The bottom map is an equivalence by the defining property of the
reflection. Since this holds for every $M$, Yoneda shows that
$\widetilde\Hsol(\eta)$ is an equivalence. Hence $\eta$ is an
$\Hsol$-equivalence.

By the preceding paragraph, $Y$ is $\Hsol$-local. Thus $\eta$
exhibits $Y$ as $\LH X$. In particular, $\LH X$ is simply connected;
and if $X$ is already $\Hsol$-local, then $\eta$ is an equivalence,
so $X$ has solid homotopy groups.
The two steps gives the equality
\[
 \Cond(\An)_{*,\mathrm{sol}}^{\ge2}
 =\Cond(\An)_{\Hsol,*}\cap\Cond(\An)_*^{\ge2}.
\]
\emph{Restricting the adjunction.}
Since $\LH$ preserves simply connected condensed anima with solid homotopy groups, the localization
adjunction restricts to pointed simply connected condensed anima.
The preceding characterization identifies its local condensed anima with
those having solid homotopy groups, proving the claimed reflection.
\end{proof}

We write $L_m$ for the left adjoint to the inclusion of the full
subcategory of pointed simply connected $m$-truncated condensed anima with solid homotopy groups
into the category of pointed simply connected condensed anima.
\begin{proposition}[Postnikov construction]\label{prop:postnikov}
For $m\ge2$,
\begin{equation}\label{eq:Postnikov}
 \tau_{\le m}\LH X\simeq L_m\tau_{\le m}X,
 \qquad \LH X\simeq\liml_m L_m\tau_{\le m}X.
\end{equation}
If $X$ is $(n-1)$-connected, $n\ge2$, then so is $\LH X$, and
$\pi_n\LH X\simeq s(\pi_nX)$.
\end{proposition}

\begin{proof}
For a
$\Hsol$-local simply connected $m$-type $Z$,
\[
 \Map_*(\tau_{\le m}\LH X,Z)\simeq\Map_*(\LH X,Z)\simeq\Map_*(X,Z)
 \simeq\Map_*(\tau_{\le m}X,Z),
\]
so $\tau_{\le m}\LH X$ reflects $\tau_{\le m}X$, which is the first formula. It also shows
that the reflected stages form a compatible Postnikov tower. Postnikov completeness then
gives the second formula, which constructs $\LH X$ from the truncated reflections alone.

If $X$ is $(n-1)$-connected, the first formula gives
$\tau_{\le n-1}\LH X\simeq L_{n-1}(*)\simeq*$ for $n\ge3$, and for $n=2$ the claim is
automatic. Thus $\LH X$ is $(n-1)$-connected. For solid $B$,
\[
\begin{aligned}
\operatorname{Hom}(\pi_n\LH X,B)
&\simeq \pi_0\Map_*(K(\pi_n\LH X,n),K(B,n))\\
&\simeq \pi_0\Map_*(\tau_{\le n}\LH X,K(B,n))\\
&\simeq \pi_0\Map_*(\LH X,K(B,n))\\
&\simeq \pi_0\Map_*(X,K(B,n))\\
&\simeq \pi_0\Map_*(\tau_{\le n}X,K(B,n))\\
&\simeq \pi_0\Map_*(K(\pi_nX,n),K(B,n))\\
&\simeq \operatorname{Hom}(\pi_nX,B).
\end{aligned}
\]
Since $\pi_n\LH X$ is solid, this identifies it with $s(\pi_nX)$.
\end{proof}

The truncated reflector $L_m$ is essential in this formula: $\LH$ itself can create homotopy
above degree $m$, so $\LH \tau_{\le m}X$ need not be an $m$-type
(Remark~\ref{rem:circle}). Only its $m$-truncation is determined by $\tau_{\le m}X$:
applying the first formula to both $X$ and $\tau_{\le m}X$ gives
$\tau_{\le m}\LH X\simeq\tau_{\le m}\LH \tau_{\le m}X$.
\\

We now turn to condensed anima that are not simply connected, where the fundamental group has
to be taken into account.
Let $G$ be a condensed group, not necessarily abelian, and
write $BG$ for its classifying condensed anima. Even if
$G^{\mathrm{ab}}$ is solid, $BG$ need not be $\Hsol$-local
(Example~\ref{ex:heisenberg}). For nilpotent condensed anima, we instead
require the fundamental group to admit a finite central series
with solid abelian factors.\\

We will use the following notion of nilpotence in a general hypercomplete topos $\mathcal X$; below it will be $\Cond(\An)$.

\begin{definition}[Nilpotent object {\cite[Definitions~A.8 and~A.10]{Mattis}}]
\label{def:nilpotent-object}
A pointed connected object $X\in\mathcal X$ is \emph{nilpotent} if $G=\pi_1X$ admits a finite central
series and, for every $n\ge2$, there is a finite filtration by $G$-stable subobjects
\[
    \pi_nX=M_0\supseteq M_1\supseteq\cdots\supseteq M_r=0
\]
such that $G$ acts trivially on each $M_i/M_{i+1}$. The length $r$ may depend on $n$.
\end{definition}
Nilpotence gives a filtration on which $G$ acts trivially on the
successive quotients. For the Postnikov argument, we need these
quotients to be solid as well.
\begin{lemma}\label{lem:invariants}
Let $G$ be a condensed group acting on a solid abelian group $M$. Then the fixed points group $M^G$ is solid. If the action is
nilpotent, then the upper invariant filtration, defined by 
\[
M^{(0)}=0 \quad \text{and} \quad  M^{(i+1)}/M^{(i)}=(M/M^{(i)})^G,
\]
is finite, $G$-stable, and exhausts $M$. Moreover, its
successive quotients are solid with trivial $G$-action.
\end{lemma}

\begin{proof}
The subgroup $M^G$ is the kernel of the map
$M\to\underline{\operatorname{Hom}}(\Z[G],M)$ sending $m$ to the function $g\mapsto gm-m$, and the target
is solid by Lemma~\ref{lem:internal-hom-solid}. Solid abelian groups are closed under kernels, cokernels
and extensions, so by induction all $M^{(i)}$ and $M/M^{(i)}$ are solid.
To show that $M^{(i)}$ is exhustive, note that if $M=M_0\supseteq\cdots\supseteq M_r=0$ is a filtration as in Definition~\ref{def:nilpotent-object}, then
induction gives $M_{r-i}\subseteq M^{(i)}$, so $M^{(r)}=M$.
\end{proof}

In the non-simply connected case one has monodromies to take into account. We therefore work with Eilenberg--Mac Lane objects
in the slice over the base.
For an abelian sheaf $M$ over a condensed anima $Y$, write
$K_Y(M,m)$ for the relative Eilenberg--Mac Lane object in the slice
over $Y$. When discussing its $\Hsol$-locality, we mean the locality of its
underlying condensed anima.

\begin{lemma}\label{lem:solid-nilpotent-gerbes}
Let $Y$ be an $\Hsol$-local condensed anima, and let $M$ be an abelian
sheaf over $Y$. Suppose that $M$ admits a finite filtration
\[
0=M_0\subset M_1\subset\cdots\subset M_r=M
\]
such that
\[
M_i/M_{i-1}\simeq Y\times A_i
\]
for solid condensed abelian groups $A_i$, where $Y\times A_i$ denotes pullback of $A_i$ from the point to $Y$. 
Then:
\begin{enumerate}
    \item $K_Y(M,m)$ is $\Hsol$-local for every $m\geq 0$.
    \item For $n\ge 1$, every $M$-banded $n$-gerbe over $Y$
    has $\Hsol$-local total space.
\end{enumerate}
\end{lemma}

\begin{proof}
We prove the first assertion by induction on $r$, simultaneously
for all $m$. For $r=0$, the relative Eilenberg--Mac Lane object is
the terminal object of the slice, whose total space is $Y$.

For the induction step we now that the claim holds for $r-1$ and prove for $r$, the short exact sequence
\[
0\longrightarrow M_{r-1}
\longrightarrow M
\longrightarrow Y\times A_r
\longrightarrow 0
\]
gives a pullback
\[
K_Y(M,m)
\simeq
\bigl(Y\times K(A_r,m)\bigr)
\times_{K_Y(M_{r-1},m+1)}Y
\]
in which the map from $Y$ to $K_Y(M_{r-1},m+1)$ is the canonical basepoint-section.  All three terms on the
right are $\Hsol$-local ($K_Y(M_{r-1},m+1)$ is local from the induction hypothesis), so their pullback is $\Hsol$-local since $\Hsol$-local condensed anima are closed under limits - the inclusion of $\Hsol$-local condensed anima is a right adjoint.

For the second assertion, apply the classification of banded
gerbes in the slice over $Y$
\cite[Theorem~7.2.2.26]{HTT}.
An $M$-banded $n$-gerbe $E\to Y$ is classified by a section
\[
\kappa:Y\longrightarrow K_Y(M,n+1),
\]
and its total space is the pullback
\[
E\simeq Y\times_{K_Y(M,n+1)}Y,
\]
formed using $\kappa$ and the base-point section. The first assertion
and closure of $\Hsol$-local condensed anima under pullbacks imply that $E$
is $\Hsol$-local.
\end{proof}
We now give a first criterion for locality beyond the simply connected case.
\begin{theorem}\label{thm:nilpotent-locality}
Let $X$ be a pointed connected condensed anima.
Assume that $B\pi_1X$ is $\Hsol$-local and that the
action of $\pi_1X$ on each $\pi_nX$, for $n\geq 2$, is nilpotent.
Then
\[
X\text{ is }\Hsol\text{-local}
\quad\Longleftrightarrow\quad
\pi_nX\text{ is solid for every }n\geq 2.
\]
\footnote{The nilpotence length may depend on $n$.}
\end{theorem}

\begin{proof}
Put $G=\pi_1X$. Suppose first that $\pi_nX$ is solid for every
$n\ge2$. By Lemma~\ref{lem:invariants}, each $\pi_nX$ admits a
finite $G$-stable filtration with solid  quotients on which $G$ acts trivially.
Starting from the $\Hsol$-local condensed anima $\tau_{\le1}X\simeq BG$, we prove
inductively that every Postnikov truncation $\tau_{\le n} X$ is $\Hsol$-local. Indeed,
$\tau_{\le n}X\to\tau_{\le n-1}X$ is a gerbe banded by the
abelian sheaf associated to the $G$-module $\pi_nX$, so
Lemma~\ref{lem:solid-nilpotent-gerbes} applies. Postnikov
completeness and closure of $\Hsol$-local condensed anima under limits imply
that $X$ is $\Hsol$-local.

Conversely, if $X$ is $\Hsol$-local, then so is its universal
cover
\[
\widetilde X=X\times_{BG}*,
\]
since $BG$ is $\Hsol$-local. By Theorem~\ref{thm:reflection},
$\pi_n\widetilde X\simeq\pi_nX$ is solid for every $n\ge2$.
\end{proof}
We next complete the above criteria and give a sufficient condition for
$B\pi_1X$ to be local.
Recall that a finite central series for $G$ is a filtration
$1=G_0\subseteq\cdots\subseteq G_c=G$ by normal subgroups such that
each $G_i/G_{i-1}$ lies in the center of $G/G_{i-1}$
\begin{corollary}\label{cor:nilpotent-solid-group}
Let $G$ be a nilpotent condensed group admitting \footnote{We require only that $G$ admit one central series with solid
factors; other central series of $G$ need not have this property.} a finite central series
\[
1=G_0\subset G_1\subset\cdots\subset G_c=G
\]
such that every quotient $G_i/G_{i-1}$ is a solid abelian group.
Then $BG$ is $\Hsol$-local.
\end{corollary}

\begin{proof}
A central extension
\[
1\longrightarrow A\longrightarrow P\longrightarrow Q
\longrightarrow 1
\]
defines a gerbe $BP\to BQ$ banded by the constant abelian
sheaf $A$. Its extension class gives a map
\[
e:BQ\longrightarrow K(A,2)
\]
and a pullback
\[
BP\simeq BQ\times_{K(A,2)}*.
\]
Thus $BQ$ $\Hsol$-local and $A$ solid imply that $BP$ is $\Hsol$-local. The result now follows by induction on the central series.
\end{proof}

Together, the preceding results give a sufficient condition
for locality of a nilpotent condensed anima. To prove necessity, we
localize a principal refinement of its Postnikov tower.
The same construction proves preservation of nilpotence.
\begin{corollary}\label{cor:nilpotent}
Let $X$ be a pointed connected nilpotent condensed anima. Then $X$ is
$\Hsol$-local if and only if $\pi_1(X)$ admits a (finite) central series with solid factors and $\pi_nX$ is solid for every $n\ge2$. Moreover $\LH X$ is
again connected and nilpotent.
    
\end{corollary}
\begin{proof}
If $\pi_1X$ admits a finite central series with solid factors,
then $B\pi_1X$ is $\Hsol$-local by
Corollary~\ref{cor:nilpotent-solid-group}.
If the higher homotopy groups are also solid,
Theorem~\ref{thm:nilpotent-locality} implies that $X$ is $\Hsol$-local.

For the converse and the final assertion, choose a principal
refinement of the Postnikov tower
\cite[Lemma~A.15]{Mattis}:
\[
X\simeq\varprojlim_i X_i,\qquad X_0=*,
\qquad
X_{i+1}\longrightarrow X_i\longrightarrow K(A_i,d_i+1),
\]
where $d_i\ge1$ and $d_i\to\infty$.
Set $Y_i=\LH X_i$ and $
Z_i=\LH K(A_i,d_i+1)
$.
Propositions~\ref{prop:fiber} give fiber sequences
\[
Y_{i+1}\longrightarrow Y_i\longrightarrow Z_i
\]
Since $\LH$ preserves connectivity, $\pi_kZ_i=0$ for
$k\le d_i$. In particular, $Z_i$ is simply connected,
and all its homotopy groups are solid.
Inductively, every $Y_i$ is connected and nilpotent, since
connected fibers of maps between nilpotent condensed anima are
nilpotent \cite[Lemma~A.12]{Mattis}.
The homotopy exact sequence shows that its higher homotopy
groups are solid. Moreover, it gives a central extension
\[
1\longrightarrow
\operatorname{coker}\bigl(\pi_2Y_i\to\pi_2Z_i\bigr)
\longrightarrow\pi_1Y_{i+1}
\longrightarrow\pi_1Y_i\longrightarrow1.
\]
Centrality follows from simple connectedness of $Z_i$.
The kernel is solid, so induction also gives a finite central
series with solid factors for $\pi_1Y_i$.

Put $Y=\varprojlim_iY_i$. 
We would like to identify $Y$ with $\LH X$.
The condensed anima $Y$ is $\Hsol$-local, being a limit of $\Hsol$-local condensed anima.
Both $X\to X_i$ and $Y\to Y_i$ become arbitrarily connected.
Ordinary chains preserve connectivity of maps, and $\Ls$
is exact and preserves connective complexes. Hence, in each
fixed homological degree and for sufficiently large $i$,
the natural comparison gives
\[
H_k\Hsol(X)
\simeq H_k\Hsol(X_i)
\simeq H_k\Hsol(Y_i)
\simeq H_k\Hsol(Y).
\]
Therefore $X\to Y$ is an $\Hsol$-equivalence, so
$Y\simeq\LH X$.
This proves preservation of connected nilpotent condensed anima;
when $X$ is $\Hsol$-local, the equivalence $X\simeq Y$ also proves
the required necessity.
\end{proof}
	If the fundamental group is abelian, the central-series condition reduces to solidity of the group itself.
\begin{corollary}\label{cor:simple}
Let $X$ be a pointed connected nilpotent condensed anima with $\pi_1X$ abelian. Then $X$ is $\Hsol$-local if and only if
$\pi_nX$ is solid for every $n\ge1$.
\end{corollary}

\section{Comparison with stable solidification}
\label{sec:stable}

One can reduces the calculation of $\LH X$ to the solidification of its
Eilenberg--Mac Lane layers. We identify these layers through the stable solidification.
Denote by $\UZ$ the forgetful functor from $D(\Cond(\Ab))$ to
$\Cond(\Sp)$,
\begin{proposition}\label{prop:stablecalc}
There is a natural equivalence
\begin{equation}\label{eq:HZ}
 \Lsp\UZ M\simeq\UZ\Ls M
 \qquad\forall M\in D(\Cond(\Ab)).
\end{equation}
\end{proposition}

\begin{proof}

Let $Y$ be a solid spectrum.
Remark that if we show that \begin{equation}\label{eq:eq44}
 \Map(\UZ\Ls M,Y)\simeq
 \Map(\UZ M,Y) \end{equation}
then since $\UZ\Ls M$ is solid, this is the defining property of $\Lsp\UZ M$ and we get \eqref{eq:HZ}.
Hence we are left to prove \eqref{eq:eq44}.  By lemma \ref{solidhom}
$\underline{\operatorname{Hom}}(H\Z,Y)$ is solid hence lies in  \[
\operatorname{Mod}_{\mathbb{Z}}(\operatorname{Solid}(\operatorname{Sp}))
\simeq
D(\operatorname{Solid}(\operatorname{Ab})).
\] The cofree-module and
derived-solidification adjunctions therefore give
\[
 \Map(\UZ\Ls M,Y)\simeq
 \Map(\Ls M,\underline{\operatorname{Hom}}(H\Z,Y))\simeq
 \Map(M,\underline{\operatorname{Hom}}(H\Z,Y))\simeq
 \Map(\UZ M,Y).
\]
This is exactly $\eqref{eq:eq44}$.
\end{proof}
To compute the solidification of Eilenberg--Mac Lane objects, we compare the unstable and stable solidification on infinite loop spaces.
\begin{proposition}[Eilenberg--Mac Lane objects]\label{prop:EM}
For $E\in\Cond(\Sp)_{\ge1}$ there is a natural equivalence
\[
\LH\Omega^\infty E\simeq\Omega^\infty\Lsp E.
\]
In particular, for $A\in\Cond(\Ab)$ and $n\ge1$,
\begin{equation}\label{eq:EM}
\LH K(A,n)\simeq\Omega^\infty\UZ\Ls(A[n]),
\qquad
\pi_{n+p}\LH K(A,n)\simeq L_p s(A)
\quad(p\ge0).
\end{equation}
\end{proposition}

\begin{proof}
Set $X_j=\Omega^\infty\Sigma^jE$ and $Y_j=\LH X_j$ for $j\ge0$.
Since $X_{j+1}$ is simply connected, Proposition~\ref{prop:fiber}
applied to $X_j\to *\to X_{j+1}$ gives
\[
Y_j\simeq\Omega Y_{j+1},
\]
compatibly with the maps $X_j \to Y_j$. For $j\ge1$, the condensed anima $Y_j$ is
$(j+1)$-connective and has solid homotopy groups. The same holds
for $j=0$ by looping $Y_1$. Thus the $Y_j$ form a $1$-connective
solid spectrum $E'$, and the maps $X_j \to Y_j$ assemble into a map of condensed spectra $E\to E'$.

Let $N$ be a solid $1$-connective spectrum. For $j\ge1$,
the condensed anima $\Omega^\infty\Sigma^jN$ is simply connected with
solid homotopy groups, hence $\Hsol$-local. Since mapping spaces between $\Omega$-spectra are computed as the limit of the maps between the components we get that
\[
\begin{aligned}
\Map(E',N)
&\simeq \liml_{j\ge1}
  \Map_*(Y_j,\Omega^\infty\Sigma^jN)\\
&\simeq \liml_{j\ge1}
  \Map_*(X_j,\Omega^\infty\Sigma^jN)\\
&\simeq \Map(E,N).
\end{aligned}
\]
For arbitrary solid $N$, the same conclusion follows by
replacing $N$ with $\tau_{\ge1}N$, since $E$ and $E'$ are
$1$-connective. Hence $E'\simeq\Lsp E$.
Taking zeroth spaces proves the first assertion.
Applying \eqref{eq:HZ} to $A[n]$ proves the remaining formulas.
\end{proof}

\begin{example}\label{ex:real}
Let $V$ be a condensed $\underline{\mathbb R}$-module. Then $\Ls V=0$
\cite[Example 2.22]{Aoki}. Indeed, the identity of $V$ factors as
$V\to\underline{\mathbb R}\otimes^{\mathbb L}V\to V$, and
$\Ls(\underline{\mathbb R}\otimes^{\mathbb L}V)\simeq
\Ls\underline{\mathbb R}\sotimes\Ls V=0$. Hence $\LH K(V,n)\simeq*$ for all $n\ge1$.
\end{example}

\begin{remark}[Failure of degreewise solidification]\label{rem:circle}
The derived solidification in \eqref{eq:EM} is necessary. For
$\mathbb T=\underline{\mathbb R/\Z}$, the exponential sequence and
$\Ls\underline{\mathbb R}=0$ \cite[Corollary 6.1(iii)]{CS} give $\Ls\mathbb T=\Z[1]$.
Consequently, for $n\ge1$,
\[
 s\mathbb T=0\qquad \text{while} \qquad \LH K(\mathbb T,n)\simeq K(\Z,n+1).
\]
Solidification can therefore move the homotopy of a single Eilenberg--Mac Lane object into
a higher degree. The case $n=1$ says that the classifying condensed anima $B\mathbb T$ localizes to
$K(\Z,2)$, so $\LH$ can kill a nontrivial $\pi_1$ and create $\pi_2$.
\end{remark}

\section{The solidification spectral sequence}
\label{sec:SS}
We now construct a spectral sequence that computes the homotopy groups of the solidification for simply connected condensed anima, by assembling the Eilenberg--Mac Lane layers along the Whitehead tower. Let
$X\in\Csc$. Its Whitehead tower has terms $\tau_{\ge q}X=\fib(X\to\tau_{\le q-1}X)$ for $q\ge2$, and. Write
\[
 F_qX=\LH(\tau_{\ge q}X),\qquad
 G_qX=\Omega^\infty\Ls((\pi_qX)[q]).
\]
The Whitehead fiber sequences $\tau_{\ge q+1}X\to \tau_{\ge q}X\to K(\pi_qX,q)$ have simply connected base.
The fiber lemma and \eqref{eq:EM} turn them into
\begin{equation}\label{eq:layers}
 F_{q+1}X\longrightarrow F_qX\longrightarrow G_qX.
\end{equation}
This is a decreasing filtration of $F_2X=\LH X$, with $F_qX$ being $(q-1)$-connected by
Proposition~\ref{prop:postnikov}. 

\begin{theorem}[Derived homotopy groups]\label{thm:SS}
There is a natural strongly convergent spectral sequence of solid abelian groups
\begin{equation}\label{eq:SS}
 E^2_{p,q}=L_p s(\pi_qX)\Longrightarrow\pi_{p+q}\LH X,
 \qquad p\ge0,\ q\ge2,
 \qquad d_r:(p,q)\mapsto(p-r,q+r-1).
\end{equation}
The induced filtration on $\pi_n\LH X$ is finite and is given by
$\im(\pi_nF_qX\to\pi_n\LH X)$, with $2\le q\le n$.
\end{theorem}

\begin{proof}
Take the homotopy exact couple of \eqref{eq:layers}:
\[
 \cdots\longrightarrow\pi_nF_{q+1}X\longrightarrow\pi_nF_qX\longrightarrow\pi_nG_qX
 \longrightarrow\pi_{n-1}F_{q+1}X\longrightarrow\cdots.
\]
In the range $n\ge2$ this consists of abelian groups, and it ends at $\pi_1F_{q+1}X=0$. It
can be read objectwise on extremally disconnected sets, where evaluation is exact. Hence it
is an exact couple of condensed abelian groups, and all of its terms are solid. The layer groups are
$\pi_{p+q}G_qX=L_p s(\pi_qX)$, and the first differential is
\[
 \pi_{p+q}G_qX\longrightarrow\pi_{p+q-1}F_{q+1}X
       \longrightarrow\pi_{p+q-1}G_{q+1}X.
\]
This has bidegree $(-2,1)$, which accounts for the initial page label $E^2$, and the derived
exact couples give the stated $d_r$. Each page is a subquotient of solid groups in the abelian
category $\Solid(\Ab)$.

In degree $n$, the filtration begins with the whole group at $q=2$ and vanishes for $q>n$,
because $F_qX$ is $(q-1)$-connected. Moreover, $d_r$ leaving $(p,q)$ vanishes for $r>p$ and $d_r$
entering $(p,q)$ vanishes for $r>q-1$, so each position stabilizes at a finite page. This gives
strong convergence, with only $2\le q\le n$ contributing to $\pi_n\LH X$.
\end{proof}

\begin{remark}
On the line $p=0$, the edge map $s(\pi_nX)\simeq\pi_nF_nX\to\pi_n\LH X$ is the factorization of
$\pi_n$ of the unit through $s(\pi_nX)$, which exists because $\pi_n\LH X$ is solid.
Thus its kernel and cokernel measure the failure of solidification to act degreewise.
\end{remark}

The spectral sequence also gives a criterion for a map to become an equivalence after solidification.
\begin{corollary}[Comparison]\label{cor:comparison}
Let $f:X\to Y$ be a map in $\Csc$ such that $\Ls(\pi_qf)$ is an equivalence for every $q\ge2$.
Then $\LH f$ is an equivalence. In particular, if $\Ls(\pi_qX)=0$ for every $q$, then $\LH X\simeq*$.
\end{corollary}

\begin{proof}
The map induces an isomorphism on $E^2$, hence on $E^\infty$. By the finiteness of the
filtrations, $\pi_n(\LH f)$ is an isomorphism for every $n$.
\end{proof}

By Example~\ref{ex:real}, for instance, $\LH X\simeq*$ whenever all the homotopy groups of $X$ are
condensed $\underline{\mathbb R}$-modules.

    The comparison with stable solidification gives a corresponding spectral sequence for bounded-below condensed spectra.
\begin{proposition}[Spectra]\label{prop:stable-SS}
For $E\in\Cond(\Sp)_{\ge a}$ there is a natural strongly convergent spectral sequence
\[
E^2_{p,q}=L_ps(\pi_qE)\Longrightarrow\pi_{p+q}\Lsp E,
\qquad p\ge0,\quad q\ge a.
\]
with the differentials of \eqref{eq:SS} and a finite filtration in each degree.
\end{proposition}

\begin{proof}
It follows from the pointed spaces version by applying $\Omega^\infty$.
\end{proof}

\section{Some   counterexamples}
\label{sec:counterexamples}
 
Let $P=\prod_{\mathbb N}\Z$. Although $P$ is solid, its ordinary condensed tensor square $P\otimes P$
is not. Indeed, solidification replaces it by the completed tensor product
(Theorem~\ref{thm:solid-tensor}):
\begin{equation}\label{eq:producttensor}
 s(P\otimes P)=P\otimes_{\square}P=\prod_{\mathbb N^2}\Z.
\end{equation}
Evaluation at a point preserves colimits and sends $\Z[T]\otimes\Z[T']=\Z[T\times T']$ to the
corresponding free group, so it computes ordinary tensor products: $(P\otimes P)(*)=P(*)\otimes P(*)$.
On these sections the unit sends a finite sum of pure tensors to a finite-rank integer matrix, so it
misses the infinite identity matrix.

\begin{example}[Solid abelianization does not imply locality]
\label{ex:heisenberg}
Let $P=\prod_{\mathbb N}\Z$ and $C=P\otimes P$.
Consider the condensed Heisenberg group
\[
H=P\times P\times C,\qquad
(a,b,c)(a',b',c')
=(a+a',\,b+b',\,c+c'+a\otimes b').
\]
It fits into a central extension
\[
1\longrightarrow C\longrightarrow H
\longrightarrow P\times P\longrightarrow1.
\]
Since
\[
[(a,0,0),(0,b,0)]=(0,0,a\otimes b),
\]
its commutator subgroup is $C$, and hence
$H^{\mathrm{ab}}\simeq P\times P$ is solid.
Nevertheless, $BH$ is not $\Hsol$-local. Indeed, $C$ is not
solid, as shown in the proof of Proposition~\ref{prop:noadjoints}.
The fiber sequence
\[
BC\longrightarrow BH\longrightarrow B(P\times P)
\]
has $\Hsol$-local base. If $BH$ were $\Hsol$-local, closure of $\Hsol$-local condensed anima
under pullbacks would imply that $BC$ is $\Hsol$-local. By
Remark~\ref{prop:EM}, this would force $C$ to be solid,
a contradiction.
\end{example}
 	The same failure of solidity for $P\otimes P$ also appears in the homotopy of a wedge.
\begin{proposition}\label{prop:noadjoints}
The inclusion
\[
 \Cond(\An)_{*,\mathrm{sol}}^{\ge2}\hookrightarrow\Cond(\An)_*^{\ge2}
\]
has no right adjoint. In fact, it does not preserve coproducts.
\end{proposition}
 
\begin{proof}
Set $W=K(P,2)\vee K(P,2)$, one can see that
\[
 \pi_3W=P\otimes P.
\]
Thus a wedge of solid simply connected condensed anima need not have solid homotopy groups. The coproduct
of the two copies of $K(P,2)$ in $\Csol$ is $\LH W$, and $W\to \LH W$ is not an equivalence. So the
inclusion does not preserve coproducts, and it has no right adjoint.
\end{proof}

    We can now compare the homotopy groups of this wedge with those of its solidification.
\begin{example}\label{ex:wedge}
As $P$ is
solid, $L_ps(P\oplus P)=0$ for $p>0$, so
\[
 \pi_2\LH W\simeq P\oplus P,\qquad
 \pi_3\LH W\simeq s(P\otimes P)\simeq P\otimes_\square P.
\]
Thus the coproduct of solid condensed anima sees the completed tensor product.
\end{example}
 
\begin{remark}
This is an unstable phenomenon: by Proposition~\ref{prop:stable}, the inclusion of solid spectra has
both adjoints.
\end{remark}
\bibliographystyle{alphaurl} 
\bibliography{meta/references} 

@book{HTT,
  title={Higher Topos Theory (AM-170)},
  author={Lurie, Jacob},
  isbn={9780691140490},
  lccn={2008038170},
  series={Annals of Mathematics Studies},
  url={https://books.google.co.il/books?id=CTe68E8wK4QC},
  year={2009},
  publisher={Princeton University Press}
}

@misc{Aoki,
  author        = {Aoki, K.},
  title         = {(Semi)topological {$K$}-theory via solidification},
  year          = {2024},
  eprint        = {2409.01462v1},
  archivePrefix = {arXiv},
  url           = {https://arxiv.org/abs/2409.01462v1}
}

@misc{BEM,
  author        = {Bachmann, T. and Engelmann, A. and Mattis, K.},
  title         = {Monadic resolutions for generalized spaces},
  year          = {2025},
  eprint        = {2511.07064v1},
  archivePrefix = {arXiv},
  url           = {https://arxiv.org/abs/2511.07064v1}
}

@misc{Mattis,
  author        = {Mattis, Klaus},
  title         = {Unstable {$p$}-completion in motivic homotopy theory},
  year          = {2024},
  eprint        = {2401.17848},
  archivePrefix = {arXiv},
  primaryClass  = {math.AG},
  url           = {https://arxiv.org/abs/2401.17848}
}

@article{Bousfield1975,
  author  = {Bousfield, A. K.},
  title   = {The localization of spaces with respect to homology},
  journal = {Topology},
  volume  = {14},
  number  = {2},
  pages   = {133--150},
  year    = {1975},
  doi     = {10.1016/0040-9383(75)90023-3}
}

@misc{CS,
  author = {Scholze, P.},
  title  = {Lectures on condensed mathematics},
  year   = {2019},
  note   = {Joint work with D. Clausen},
  url    = {https://www.math.uni-bonn.de/people/scholze/Condensed.pdf}
}

\end{document}